\documentclass[12pt, oneside]{article}

\usepackage[dvips]{graphicx}
\usepackage{amsthm}
\usepackage{amsmath}
\usepackage{amsfonts}
\usepackage{amssymb}

\theoremstyle{plain}
\newtheorem{thm}{\textbf{Theorem}}[section]

\newtheorem{con}{\textbf{Conjecture}}[section]

\newcommand{\ignore}[1]{}
\begin{document}
\title{Domination Cover Pebbling: Structural Results}

\author{
Nathaniel G. Watson\\
Department of Mathematics\\
Washington University at St.~Louis\and
Carl R.~Yerger\\
Department of Mathematics\\
Georgia Institute of Technology}
\maketitle

\begin{abstract}This paper continues the results of ``Domination Cover
Pebbling: Graph Families.'' An almost sharp bound for the
domination cover pebbling (DCP) number, $\psi(G)$, for graphs $G$
with specified diameter has been computed.  For graphs of diameter
two, a bound for the ratio between $\lambda(G)$, the cover
pebbling number of $G$, and $\psi(G)$ has been computed.  A
variant of domination cover pebbling, called subversion DCP is
introduced, and preliminary results are discussed.
\end{abstract}

\section{Introduction}
\phantom{space } Given a graph $G$ we distribute a finite number of indistinguishable markers
 called \emph{pebbles} on its vertices. Such an arrangement of pebbles, which can also be thought of as a function from $V(G)$ to
$\mathbb{N} \cup \{0\},$ is called a
\textit{configuration}. A \emph{pebbling move} on a graph is defined as taking two pebbles off one vertex, throwing one away, and moving the other to an adjacent vertex. Most research in pebbling has focused on a quantity known as the \emph{pebbling number} $\pi(G)$ of a graph, introduced by F. Chung in \cite{Chung}, which is defined to be the smallest integer $n$ such that for every configuration of $n$ pebbles on the graph and for any vertex $v \in G,$ there exists a sequence of pebbling moves starting at this configuration and ending in a configuration in which there is at least one pebble on $v$. A new variant of this concept, introduced in by Crull et al.\ in \cite{Crull}, is the \emph{cover pebbling number} $\lambda(G)$, defined as the minimum number $m$ such that for any initial configuration of at least $m$ pebbles on $G$ it is possible to make a sequence of pebbling moves after which there is at least one pebble on every vertex of $G$.

In a recent paper (\cite{VNIDCP1}) the authors, along with
Gardner, Godbole, Teguia, and Vuong, have introduced a concept
called domination cover pebbling and have presented some
preliminary results.
Given a graph $G,$ and a configuration $c,$ we call a vertex $v \in G$ \emph{dominated} if it is covered (occupied by a pebble) or adjacent to a covered vertex. We call a configuration $c'$ \emph{domination cover pebbling solvable}, or simply \emph{solvable}, if there is a sequence of pebbling moves starting at $c'$ after which every vertex of $G$ is dominated.  We define the \emph{domination cover pebbling number} $\psi(G)$ to be the minimum number $n$ such that any
initial configuration of $n$ pebbles on $G$ is domination cover pebbling solvable.



The set of covered vertices in the final configuration
depends, in general, on the initial configuration---in particular,
$S$ need not equal a minimum dominating set. For instance, consider
the configurations of pebbles on $P_4$, the path on four vertices,
as shown in Figure 1:
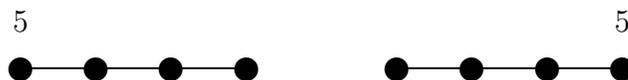
\begin{figure}[htb] \unitlength 1mm
\begin{center}
\begin{picture}(100,30)
\put(5,5){\circle*{3}} \put(15,5){\circle*{3}}
\put(25,5){\circle*{3}}
\put(35,5){\circle*{3}}\put(55,5){\circle*{3}}
\put(65,5){\circle*{3}}\put(75,5){\circle*{3}}
\put(85,5){\circle*{3}}
\put(5,5){\line(1,0){30}} 
\put(55,5){\line(1,0){30}} \put(4,10){5} \put(84,10){5}
\end{picture}
\end{center}
\caption{An example where two different initial configurations
produce two different domination cover solutions.} \label{ex2}
\end{figure}

For the graph on the left, we make pebbling moves so that the first
and third vertices (from left to right) form the vertices of the
dominating set.  However, for the graph on the right, we make
pebbling moves so that the second and fourth vertices are selected
to be the vertices of the dominating set.  In some cases, moreover,
it takes more vertices than are in the minimum dominating set of
vertices to form the domination cover solution.  For example, in
Figure 2 we consider the case of the binary tree with height two,
where the minimum dominating set has two vertices, but the minimal
dominating set possible for a domination cover solution has three
vertices. This corresponds to several possible starting
configurations, for example the configuration pictured, the configuration with a
pebble at the leftmost bottom vertex and 4 pebbles at the
root, and the configuration with 1 and 10 pebbles at the leftmost and rightmost
bottom level vertices respectively.
\begin{figure}[tree]
\unitlength 1mm
\begin{center}
\begin{picture}(30,25)
\put(15,25){\circle*{3}} \put(15,25){\line(-1,-1){11}}
\put(15,25){\line(1,-1){11}}
\multiput(5,15)(20,0){2}{\circle*{3}}
\put(5,15){\line(-1,-2){5}}
\put(5,15){\line(1,-2){5}} \put(25,15){\line(-1,-2){5}}
\put(25,15){\line(1,-2){5}}
\multiput(0,5)(10,0){2}{\circle*{3}}
\multiput(20,5)(10,0){2}{\circle*{3}}
\put(-.8,-.5){1} \put(4,9){1} \put(24.3,9){1}
\end{picture}
\end{center}
\caption{A reachable minimal configuration of pebbles on $B_2$
that forces a domination cover solution.} \label{ex1}
\end{figure}
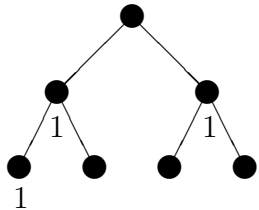

The above two facts constitute the main reason why domination cover
pebbling is nontrivial.  We refer the reader to \cite{haynes} for
additional exposition on domination in graphs, and to \cite{VNIDCP1} for some further explanation of the domination cover pebbling number, including the computation of the domination cover pebbling number for some families of graphs.


One way to understand the size of the numbers $\pi(G), \lambda(G),$ and $\psi(G)$
is to find a bound for the size of these numbers given the diameter of
$G$ and the number of vertices. This has been done for $\pi(G)$ for graphs of diameter two in \cite{Clarke} and for graphs of diameter three in $\cite{Bukh}.$ A theorem proven in \cite{jonas} and \cite{stacking} gives as a corollary a sharp bound for graphs of all diameters, which was originally established by other means in \cite{firstpaper}. In this paper, we prove that for
graphs of diameter two with $n$ vertices, $\psi(G) \leq n-1$. For
graphs of diameter $d,$ we show $\psi(G) \leq 2^{d-2}(n-2)+1$. We
also compute that the ratio $\lambda(G) / \psi(G) \geq 3$ for graphs
of diameter two.

Another way to extend cover pebbling is called subversion
domination cover pebbling.  A parameter $\omega$ used in
calculating the vertex neighbor integrity of a graph $G$ counts
the size of the largest undominated connected subset of $G$.  When
$\omega = 0$, this corresponds to domination cover pebbling.  To
conclude this paper, we provide some preliminary results for this
generalized parameter.

\section{Diameter Two Graphs}
In the next few sections, we will present structural domination cover
pebbling results.
\begin{thm} \label{dia2}
For all graphs $G$ of order $n$ with maximum diameter two, $\psi(G) \leq
n-1$.
\end{thm}

\begin{proof}
First, we show this bound is sharp by exhibiting a graph $G$ such
that $\psi(G)
> n-2$. Consider the star graph on $n$ vertices, and place a pebble
on all of the outer vertices except one.  This configuration of
pebbles does not dominate the last outer vertex.  Hence, $\psi(G) >
n-2$.

To prove the theorem, we will show that, given a graph $G$ of diameter two on $n$ vertices, any configuration $c$
of $n-1$ pebbles on $G$ is solvable.

Given such a graph configuration $c$, let
$S_1$ be the set of vertices $v \in G$ such that $c(v)
> 1$. Let $S_2$ be the set vertices $w \in G$ such that $c(w)= 0$ and $w$ is adjacent to some vertex of $S_1,$ and let $S_3$ be the rest of the vertices, the ones that are neither in $S_1$ nor adjacent to a vertex of $S_1$.
Let $a := |S_2|$, and $b := |S_3|$.
Given a configuration $c'$, define the \emph{pairing number} $P(c')$ to be  $\sum_{v \in G} $ $\max{\{0, \frac{c'(v)-1}{2}\}}$.
It can easily be checked
that $P(c')= \frac{a+b-1}{2}.$ Note that if $P(c') = k$ then $c'$
contains at least $\lceil k \rceil$ disjoint pairs of pebbles, which means
that we can make at least $\lceil k \rceil$ pebbling moves. Also, note that every vertex in $G$ is at distance at most two from some vertex in $S_1.$ This ensures that that every vertex in $S_3$ is adjacent to a vertex in $S_2.$ Also, if some vertex in $S_1$ is not adjacent to a vertex of $S_2$, it must be adjacent only to vertices in $S_1$. Since this vertex has distance at most two from any other vertex on the graph, we conclude that every vertex of the graph is either in $S_1$ or adjacent to a vertex of $S_1,$ meaning the $G$ is already dominated by covered vertices, as desired. Therefore, it suffices to consider the case in which $S_2$ is a dominating set of $G.$

First, suppose that $a \leq b$.  In this case, $P(c) \geq \frac{2a-1}{2}$.
Hence, there are at least $a$ disjoint pairs of pebbles that can be moved
from elements in $S_1$ to $S_2$. For each uncovered vertex $v \in S_2$, if
possible, move a pair of pebbles from an adjacent element of $S_1$ to put a pebble on $v$. After this is done for as many vertices of $S_2$ as possible, let $L$ be the set vertices in $S_2$ which are still uncovered. Note that these vertices are necessarily at distance $2$ from all remaining pairs of pebbles. Furthermore, since $S_1$ initially had at least $a$ disjoint pairs of pebbles, there remain at least as many pairs as there are vertices in $L.$ If this number is $0,$ the dominating set $S_2$ is covered and we are done. Otherwise, we nonetheless now
know $S_3$ is dominated because if there were some vertex $y$ that
were adjacent to only those elements of $S_2$ which are also in $L,$
then the minimum distance between $y$ and a vertex in $S_1$ with a pair of
pebbles is $3$, which is
impossible. However, it may be the case for some $z \in L$ that the
vertex in $S_1$ that $z$ was adjacent to lost its pebbles, and if
this is the case, move a pair of pebbles from $S_1$ so that $z$ is
dominated (this always possible since our graph has diameter two).
With the $|L|$ pairs we of pebbles we have, we can ensure each
vertex of $L$ is dominated. After this is done, $G$ will be
completely dominated by covered vertices.

Now consider the case $a > b$. We know that
$P(c) \geq \frac{2b-1}{2}$ and so there are at least
$b$ pairs of pebbles available. Given any vertex $v$ in $S_3$ and a
pair of pebbles on a vertex $w \in S_1$, we can use this pair to
move to a vertex between $v$ and $w,$ which is clearly in $S_2.$ We
now do this whenever necessary for each vertex of $S_3,$ first using
those pairs which can be removed from vertices having at least 3
pebbles. Let $m$ be the number of moves that have been made. Then we
know that $m$ vertices in $S_2$ now have pebbles on them.
Furthermore we know $m \leq b,$ and since some of our moves may
dominate multiple vertices of $S_3,$ thus making some other moves
unnecessary, it is indeed possible that $m<b.$ In any case, after
the moves are made, every vertex in $S_3 \cup S_1$ is dominated. If
every vertex we have removed pebbles from is still covered, then the
vertices of $S_2$ are still dominated and we are done.

Otherwise, we have removed pebbles from some vertex which had
exactly two pebbles on it. Thus, these first $m$ pebbling moves
subtract at most $\frac{2m-1}{2}$ from $P(c)$, leaving a pairing
number of $\frac{a+b-2m}{2}\geq \frac{a-m}{2}$ for the configuration
after these moves. At this point, since we were forced to use
pebbles from a vertex that had only two pebbles, we know that every
vertex that contributes to the pairing number has exactly two
pebbles on it. Thus there are at least $a-m$ vertices in $S_1$ with
two pebbles on them. We can use these pairs to dominate the $a-m$
vertices of $S_2$ which are not covered. This leaves $G$ dominated
by covered vertices and therefore $\psi(G) \leq n-1$.
\end{proof}

We can apply this theorem to prove a result about the ratio between
the cover pebbling number and the domination cover pebbling number
of a graph. We conjecture that this ratio holds for all graphs, but
it does not seem that this can be directly proven using the
structural bounds in this paper.
\begin{thm}
For all graphs $G$ of order $n$ with diameter two, $\lambda(G) / \psi(G)
\geq 3$.
\end{thm}
\begin{proof}
First, suppose that the minimum degree of a vertex of $G$ is less
than or equal to $\lceil \frac{n-1}{2} \rceil$.  By the previous
theorem, we know that the maximum value of $\psi(G)$ is $n-1$. We
now construct a configuration of pebbles on $G$ such that
$\lambda(G) \geq 3n - 3$.  Place $3n-3$ pebbles on any vertex $v$
that has a degree less than $\lceil \frac{n-1}{2} \rceil$.  It
takes $2$ pebbles to cover solve each vertex adjacent to $v$, at
most $\lceil \frac{n-1}{2} \rceil$, and all the remaining vertices
require $4$ pebbles. Since there are at least as many vertices a
distance of $2$ away from $v$ as there are a distance of $1$ away
from $G$, $3n-3$ pebbles or more are required to cover pebble all
of the vertices except for $v$.  Thus for this class of graphs,
$\lambda(G)
> 3n -3 \geq 3 \psi(G)$.

Now suppose that the minimum degree $k$ of a vertex in $G$ is
greater than $\lceil \frac{n-1}{2} \rceil$.  By a similar argument
as the previous paragraph, notice that $\lambda(G)$ for any
diameter two graph is at least $4n - 2m - 3$, where $m$ is the
minimum degree of a vertex of $G$. Since
$\lambda(G) \geq 4n - 2m - 3$, it suffices to show we can always solve a configuration $c$ of $\lfloor \frac{4n
- 2m - 3}{3} \rfloor = k$ pebbles on $G$. Given a
particular value for $m$ between $\lceil \frac{n + 1}{2} \rceil$
and $n -1$, we will construct a domination cover solution.

As long as there exist vertices of $G$ that have at least three pebbles and are adjacent to an unoccupied vertex, we haphazardly make moves from such vertices to adjacent unoccupied vertices. We claim that the resulting configuration has the desired property that the set of occupied vertices are a dominating set of $G$. First suppose that the algorithm is forced to terminate while there remains some vertex $v$ having at least three pebbles. Then this vertex must be adjacent only to occupied vertices of $G,$ and since the diameter of $G$ is two, these neighbors $v$ form a dominating set of $G$. Otherwise, if every vertex has less than three pebbles, it can easily be checked that the number of occupied vertices is now  $\sum_{v \in G} \lceil \frac{c(v)}{2} \rceil \geq \lceil \frac{k}{2} \rceil.$ Since the minimum degree of a vertex in $G$ is $m,$ by the pigeonhole principle, if we now have $n - m$ or more vertices covered by a pebble, then every vertex of $G$ is dominated.  So if $\lceil \frac{k}{2} \rceil \geq n-m$, we are
finished.  We see that
\begin{equation*} \left\lceil \frac{\left\lfloor \frac{4n - 2m - 3}{3} \right\rfloor}{2}
\right\rceil \geq \left\lceil \frac{ \frac{4n - 2m - 5}{3} }{2}
\right\rceil = \left\lceil \frac{4n}{6} - \frac{m}{3} -
\frac{5}{6} \right\rceil
\end{equation*} Therefore, we are done if

\begin{equation*} \left\lceil \frac{4n}{6} - \frac{m}{3} - \frac{5}{6}  \right\rceil \geq
n-m, \end{equation*}  which is equivalent to
\begin{equation*} n \leq \left\lceil \frac{4n}{6} + \frac{2m}{3} -
\frac{5}{6} \right\rceil.\end{equation*} This inequality holds for
$m \geq \lceil \frac{n + 1}{2}\rceil$. Therefore, we have
completed this case and have shown that for all graphs $G$ of
diameter two, $\lambda(G) / \psi(G) \geq 3$.
\end{proof}
\noindent We now prove a more general bound for graphs of diameter
$d$.
\section{Graphs of Diameter $d$}

\begin{thm}

Let $G$ be a graph of diameter $d \geq 3$ and order $n$. Then $\psi(G) \leq 2^{d-2}(n-2)+1.$
\end{thm}
 Throughout the proof, we adopt the convention that if $G$ is a graph and $V$
and $W$ are subsets of $V(G)$ and $v \in V(G)$ then $d(v,W)= \min_{w
\in W} d(v,w)$ and $d(V,W)= \min_{v \in V} d(v,W).$ Also, for any
set $S \subseteq V(G)$ we of course let $S^C = V(G) \setminus S.$
\begin{proof}
First, we define the \emph{clumping number}
$\chi$ of a configuration $c'$ by $$\chi(c') := \sum_{v \in \, G}
2^{d-2}\max\left( \left\lfloor \frac{c'(v)-1}{2^{d-2}}\right\rfloor, \ 0
\right).$$ The clumping number counts the number of pebbles in a
configuration which are part of disjoint ``clumps'' of size $2^{d-2}$ on a
single vertex, with one pebble on each occupied vertex ignored.

Now let $c$ be a configuration on $G$ of size at least
$2^{d-2}(n-2)+1.$ We will show that $c$ is solvable
 by giving a recursively defined algorithm for solving $c$ through a sequence of
pebbling moves. First, we make some definitions to begin the
algorithm:
\begin{itemize}
\item $c_0=c$.
\item $A_0 = \{ v  \in G \ : \ c(v) > 0  \}$.
\item $B_0 = \{ v \in G \ : \ c(v) \geq 2^{d-2}+1 \}$.
\item $C_0 = V(G) - A_0$.
\item $D_0 = \emptyset$.
\end{itemize}

We will describe our algorithm by recursively defining a sequence of configurations $c_p$  and four sequences $ A_p, B_p, C_p, $ and $D_p$ of sets of vertices. At each step, we will need to make sure a few
conditions hold, to ensure that the next step of the algorithm may
be performed.  For each $m$, we will insist that:
\begin{enumerate}

\item For every $v \in C_m \cup D_m$, $c_m(v) = 0$ and for every $v \in A_m,$ $c_m(v) >
0$.
\item $\chi(c_m) \geq 2^{d-2}(|C_m| - 1)$.
\item $|C_m| \leq |C_0| - m$.
\item $B_m = \{ v \in G \ : \ c_m(v) \geq 2^{d-2}+1 \}$.
\item If both $B_m \not= \emptyset$ and $D_m \not= \emptyset,$ $d(B_m, D_m) = d$ ; If $D_m \not= \emptyset,$ there always exists some $v \in G$ such that $d(v, D_m) = d,$ even if $B_m =
\emptyset$.
\item $A_m, C_m,$ and $D_m$ are pairwise disjoint and $A_m \cup C_m \cup D_m =
V(G)$.
\item Every vertex of $D_m$ is dominated by $c_m$.
\item There exists a sequence of pebbling moves transforming $c$ to $c_m$.
\end{enumerate}
Note by 1, 4, and 6, we will always have $B_m \subseteq A_m.$ Also, by 1,
6, and 7, every vertex of $G$ which is not dominated by $c_m$ is in $C_m.$

For $m=0$, only condition 2 is not immediately clear. To verify it, note
that \begin{eqnarray*} \chi(c) &=& \sum_{v \in G} 2^{d-2} \max \left(
\left\lfloor \frac{c(v)-1}{2^{d-2}}\right\rfloor, \ 0 \right) \\ & =&
\sum_{v \in A_0} 2^{d-2} \left\lfloor \frac{c(v)-1}{2^{d-2}}\right\rfloor
\\
& \geq & \sum_{v \in A_0}
2^{d-2}\left(\frac{c(v)}{2^{d-2}}-1\right).
\end{eqnarray*} Using the fact that the size of $c$ is at least
$2^{d-2}(n-2)+1,$ and $|C_0|=n-|A_0|,$ we see
$$\chi(c) \geq (2^{d-2}(n-2)+1)-2^{d-2} |A_0|=2^{d-2}(|C_0|-2)+1.$$
From the definition of $\chi$, it is apparent that $2^{d-2} | \,
\chi(c).$ Thus, we indeed must have $$\chi(c)=\chi(c_0) \geq
2^{d-2}(|C_0|-1).$$

Suppose for some $p-1 > 0$ we have defined $c_{p-1}, A_{p-1},
B_{p-1}, C_{p-1},$ and $D_{p-1}$ and the above conditions hold when
$m=p-1$. We shall assume that there is some vertex in $C_{p-1}$
which is not dominated by $c_{p-1},$ for otherwise, by conditions 6,
7 and 8, $c$ is solvable and we are done. Thus $|C_{p-1}| \geq 1.$
But suppose $|C_{p-1}| = 1.$ Call this single vertex $v.$ Since it
is non-dominated, it is adjacent to only uncovered vertices. These
vertices cannot be in $C_{p-1}$ for $|C_{p-1}|=1,$ and they are not
in $A_{p-1},$ because every vertex in $A_{p-1}$ is covered by
property 1. So every vertex adjacent to $v$ is in $D_{p-1}.$ Invoke
property 5 to choose a $w \in G$ for which $d(w, D_{p-1}) = d.$ Any
path from $w$ to $v$ passes through one of the vertices in $D_{p-1}$
which is adjacent to $v,$ and is thus of length at least $d+1,$ so
$d(w, v) \geq d+1,$ contradicting the assumption that $G$ has
diameter $d.$ We have now shown that, if $C_{p-1}$ has a
non-dominated vertex, then $|C_{p-1}| \geq 2.$ In this case, we will
have $\chi(c_{p-1}) \geq 2^{d-2},$ ensuring the existence of some
clump of size $2^{d-2},$ and thus that $B_{p-1}$ is non-empty.
Therefore, we will always implicitly assume that $B_{p-1} \not=
\emptyset$. \newline

\textbf{Case 1:} $d(B_{p-1}, C_{p-1}) \leq d-2$

In this case, we choose $v' \in B_{p-1}$ and $w' \in C_{p-1}$ for
which $d(v',w') \leq d-2$  and move $2^{d(v',w')}$ pebbles from $v'$
to $w',$ leaving one pebble on $w'$ and at least one on $v'.$ We let
$c_p$ be the configuration of pebbles resulting from this move. Let
$C_p = C_{p-1} \setminus w'.$  Thus $|C_p| = |C_{p-1}| - 1 \leq
|C_0|-(p-1)-1$ and we see that condition 3 holds when $m=p.$
Furthermore, We have used at most one clump of $2^{d-2}$ pebbles so
$$\chi(c_p) \geq \chi(c_{p-1}) - 2^{d-2} \geq 2^{d-2}(|C_{p-1}|-1)
-2^{d-2} = 2^{d-2}(|C_p|-1)$$ and therefore condition 2 holds for
$p.$ Also, we let $A_p = A_{p-1} \cup \{w'\},$ let $C_p=C_{p-1} \
w',$ and $D_p=D_{p-1}$ (now, clearly condition 6 holds.) We again
let $B_p = \{ v \in G \ : \ c_p(v) \geq 2^{d-2}+1 \},$ which simply
means that we have possible removed $v'$ from $B_{p-1}$ if $v'$ now
has less than $2^{d-2}+1$ pebbles. Thus $B_{p} \subseteq B_{p-1},$
and now 1, 4, 5, 7, and, 8 are all easily seen to hold for $m=p.$
\newline

\textbf{Case 2:} $d(B_{p-1}, C_{p-1}) \geq d-1.$

If every vertex in $C_{p-1}$ is dominated by $A_{p-1},$ we are done.
Otherwise, let $w'$ be some non-dominated vertex in $C_{p-1}.$ Clearly,
 $w'$ is at distance $d-1$ or $d$ from $B_{p-1}.$ Suppose
$d(B_{p-1}, w') = d-1.$ Then $w'$ is adjacent to some (non-covered)

vertex $w''$ at distance $d-2$ from $B_{p-1}.$ By condition 1, every vertex of
$G$ which is not covered by $c_{p-1}$ is in
$C_{p-1} \cup D_{p-1}.$ But $d(B_{p-1}, C_{p-1}) \geq d-1$ and by 5,
$d(B_{p-1}, D_{p-1}) = d$ so $w'' \notin C_{p-1} \cup D_{p-1}.$ This contradiction means that $d(w',B_{p-1}) \not= d-1$
and so $d(w',B_{p-1}) = d.$

Choose some vertex in $B_{p-1}$ and call
it $v'.$ We know $d(v',w')=d$ so consider some path of length $d$
from $v'$ to $w'.$ Let $v^*$ be the unique point on this path for
which $d(v^*,v' = d-2).$ Thus $v^* \notin C_{p-1} \cup D_{p-1}$ and so
$v^* \in A_{p-1},$ and also $d(v^*, w') =2.$ Let $w''$ be some
vertex which is adjacent to both $v^*$ and $w'$ so that
$d(v',w'')=d-1.$ Then because $w''$ is uncovered (else $w'$ would be
dominated), it must be in $C_{p-1}.$ This  also means that $v^* \notin
B_{p-1}$ by the assumption that $d(B_{p-1}, C_{n-1}) \geq d-1.$

We now move one clump of $2^{d-2}$ pebbles from $v'$ to $v^*,$
adding one pebble to $v^*,$ which now, by condition 1, has at least
two pebbles. We then move two pebbles from $v^*$ and cover $w''$
with one pebble. We let $c_p$ be the configuration resulting from
these moves. We let $D_p= D_{p-1} \cup \{ w' \}$ and we again let
$B_p = \{ v \in G \ : \ c_p(v) \geq 2^{d-2}+1 \},$ which just means
we have possibly removed $v'$ from $B_{p-1},$ so $B_p \subseteq
B_{p-1}.$ If now $c_p(v^*) = 0,$ we let $A_p = A_{p-1} \cup \{w''\}
\setminus v^* \}$ and $C_p = C_{p-1} \cup \{v^*\} \setminus \{w',
w''\}.$ Otherwise, if $c_p(v^*) > 0$, let $A_p =
A_{p-1} \cup \{w''\}$ and $C_p = C_{p-1} \setminus \{w', w''\}.$
This ensures that conditions 1 and 6 still hold for $m=p.$ Also, $|C_p| \leq
|C_{p-1}|-1 \leq |C_0|-(p-1)-1$ and so condition 3 holds for $m=p.$
Furthermore, we have used only one clump of $2^{d-2}$ pebbles, because $v^*
\notin B_{p-1}$ and so by using a pebble from $v^*,$ we could not have
destroyed a clump. Thus  $$\chi(c_p) = \chi(c_{p-1}) - 2^{d-2} \geq
2^{d-2}(|C_{p-1}|-1) -2^{d-2} \geq 2^{d-2}(|C_p|-1) $$ and therefore
condition 2 holds for $p.$ Condition 5 also still holds for $m=p$
because $B_p \subseteq B_{p-1}$ and because we have added only the vertex $w'$ to $D_{p-1}$ and
$d(B_{p-1}, w') =d,$ so $d(B_{p-1}, D_p) = d.$  To see condition 7 is
still true, note that to get $D_p$ we have only added $w'$ to
$D_{p-1},$ and certainly, $w'$ is adjacent to $w'',$ which is
covered by $c_p,$ so $w'$ is dominated by $c_p.$ Also, the only
previously covered vertex of $G$ which is now uncovered is
(possibly) $v^*$ but $d(v^*, B_{p-1} )= d-2,$ and so $v^*$ is not
adjacent to any vertex in $D_{p-1}$ for, by 5, $d(B_{p-1}, D_{p-1})
= d.$ Thus, by possibly uncovering $v^*,$ we did not cause any
vertex in $D_{p-1}$ to become undominated, so 7 still holds for
$m=p$. Finally, the fact that conditions 4 and 8 still hold for $m=p$ is easily seen.
\newline

The algorithm continues as long as there is some non-dominated
vertex in $C_p.$ By condition 3, it must terminate after at most
$|C_0|$ steps, with $ |C_k| = 0$ for some $k \leq |C_0|.$ The
configuration $c_k$ clearly dominates every vertex of $G, $ and by
property 8, $c_k$ is reachable from $c$ by pebbling moves, so $c$ is
solvable.
\end{proof}

For $d\geq 3,$ Figure 3 shows a graph $G$ which is an example of a
graph of diameter $d$ with $n = 2m+d-2$ vertices for which
$\psi(G)$ comes close to the upper bound of $2^{d-2}(n-2)+1 =
2^{d-1}m + 2^{d-2}(d-2)+1.$
\begin{figure}[htb]
\unitlength 1mm
\begin{center}
\begin{picture}(70,83)


 \multiput(10,65)(10,0){4}{\circle*{2}}
 \multiput(10,80)(10,0){4}{\circle*{2}}

\put(60,80){\circle*{2}}
 \put(44.6,65){$\ldots$}
 \put(44.6,80){$\ldots$}
\put(60,65){\circle*{2}} \put(3,62.5){\framebox(63,5)}


\multiput(10,65)(10,0){4}{\line(0,1){15}}
\put(60,65){\line(0,1){15}}

\put(35,50){\line(1,3){5}} \put(35,50){\line(-1,3){5}}
\put(35,50){\line(-1,1){15}} \put(35,50){\line(-5,3){24}}
\put(35,50){\line(5,3){25}}

\multiput(35,50)(0,-15){4}{\circle*{2}}

\multiput(35,50)(0,-15){2}{\line(0,-1){15}}

\put(35,20){\line(0,-1){4}} \put(34.5, 10.65){$\vdots$}
\put(35,5){\line(0,1){4}}
\put(29,49){$u_1$}
 \put(29,34){$u_2$}

\put(29,19){$u_3$} \put(25.5,4){$u_{d-1}$}

\put(4, 64){$v_1$} \put(14,64){$v_2$}\put(24, 64){$v_3$}
\put(34,64){$v_4$}  \put(52.3,64){$v_m$}

\put(3.5, 79){$w_1$} \put(13.5,79){$w_2$}\put(23.5, 79){$w_3$}
\put(33.5,79){$w_4$}  \put(52.3,79){$w_m$}

\end{picture}
\end{center}
\caption{\label{badgraph}A graph with high DCP number. The box
represents the fact that there is an edge between every pair of
vertices inside, making the subgraph induced by $\{v_1, v_2, \ldots,
v_m \}$ a complete graph on $m$ vertices.   }
\end{figure}
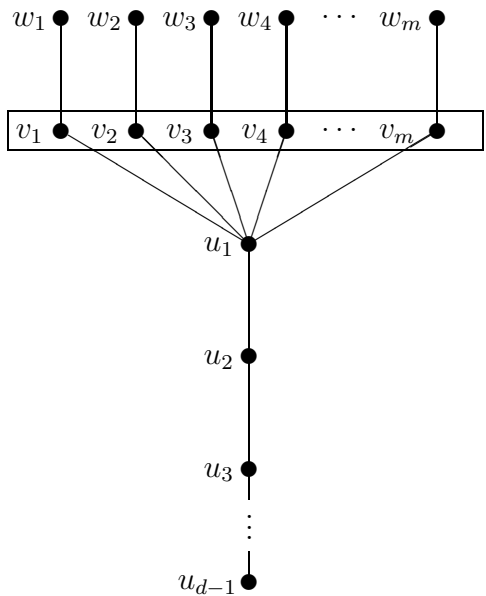

To dominate vertex $w_i,$ it is easy
to see a pebble is needed on $w_i$ or $v_i.$ They each have
distance not less than $d-1$ from $u_{d-1},$ and so it requires
$2^{d-1}$ pebbles on $u_{d-1}$ to supply this pebble. This means
at least $2^{d-1}m$ pebbles are needed on $u_{d-1}$ to dominate
every $w_i,$ so $\psi(G) \geq 2^{d-1}m.$ Further, using the result
of \cite{jonas} and \cite{stacking}, we can calculate $\lambda(G)
= 3 \cdot 2^{d-1}m+ 2^d-1.$ Clearly, by making $m$ large we can
make $\lambda(G) / \psi(G)$ arbitrarily close to 3. Also note that
for the complete graph on 2 vertices, $\lambda(G)=3$ and
$\psi(G)=1.$ We conjecture that it is not possible, however, for
the ratio to be less than 3:

\begin{con}
$\lambda(G) / \psi(G) \geq 3$ for all graphs $G$ with more than
one vertex.
\end{con}

\section{Subversion DCP}

There are several possible generalizations of domination cover pebbling which readily suggest themselves, and many of these are indeed
interesting. For instance, we may ask what happens if we simply allow $n$ vertices to remain undominated, that is, if we say a graph has been solved if all but $n$ vertices are dominated by covered vertices. More interestingly, one may relax the requirement that a graph must be dominated by pebbled vertices in order to be solved to the
condition that every vertex of a solved graph must have distance no more than $n$ from some pebbled vertex. On the other hand, we could tighten the condition that every vertex of a solved graph is either covered by pebbles or adjacent to a covered vertex by insisting that all vertices, covered or not, must be adjacent to some covered vertex.

However, these generalizations, while natural, may not be different enough from DCP to warrant extensive study. For instance, the problem of diameter bounds seems highly likely to be solvable in each case by an approach quite similar to that in Section 3. Furthermore, in each case, lower bounds which intuitively seem good can be derived from graphs quite similar to the one shown in Figure 3. Therefore, we introduce in this section a less obvious generalization of DCP which we feel makes the analogues to the questions answered in this paper more interesting than they are for the generalizations named above.


%
Given a graph $G$ and a subset $S \subseteq V(G)$, call the subgraph induced by the set of vertices which are neither in $S$ nor adjacent to a vertex of $S$ the $\emph{undominated subgraph}$ of $S$. Then we let the \emph{ $\omega$-subversion number} of $G,$ denoted $\Omega_{\omega}(G),$ be the minimum number of pebbles required
such that regardless of their initial configuration it is always possible through a sequence of pebbling moves to cover some subset of $G$ that has an undominated subgraph in which there is no connected component of more than $\omega$ vertices.\footnote{This definition and the term ``subversion" are partly inspired by Cozzens and Wu \cite{shu-shih}. Specifically, our parameter $\omega$ matches with their use of $\omega$ for the order of the largest connected component of an undominated subgraph. }  Notice
that domination cover pebbling corresponds to the case when
$\omega = 0$.

\section{Basic Results}

\begin{thm}
For $\omega \geq 0$, $\Omega_{\omega}(K_n) = 1$.
\end{thm}
\begin{proof}
When any pebble is placed on $K_n$, the entire graph is dominated.
\end{proof}
\begin{thm}
For $ s_1 \geq s_2 \geq \cdots \geq s_r$, let $K_{s_1, s_2, \ldots, s_r}$
be the complete $r$-partite graph with $s_1,s_2,\ldots,s_r$ vertices in
vertex classes $c_1, c_2, \ldots, c_r$
 respectively.  Then for $\omega \geq 1,$ $\Omega_{\omega}(K_{s_1, s_2, \ldots, s_r}) = 1$.
\end{thm}
\begin{proof}
Place a pebble on any vertex in $c_i$.  All the vertices in the
other $c_i$'s are dominated.  The other vertices in $c_1$ that are
undominated are disjoint from each other.  Thus, the result follows.
\end{proof}

\begin{thm}
For $\omega \geq 1$, $n \geq \omega + 3$, $\Omega_{\omega}(W_n) =
n-2-\omega$, where $W_n$ denotes the wheel graph on $n$ vertices.
\end{thm}
\begin{proof}
First, we will show that $\Omega_{\omega}(W_n) > n-3-\omega$.  Place
a single pebble on each of $n-3-\omega$ consecutive outer vertices
so that all of the pebbled vertices form a path.  This leaves a
connected undominated set of size $\omega + 1$. Hence,
$\Omega_{\omega}(W_n) > n-3-\omega$.  Now, suppose that we place
$n-2-\omega$ pebbles on $W_n$.  If any vertices have a pair of
pebbles on them, the entire graph can be dominated by moving a
single pebble to the hub vertex.  Hence, each vertex can contain
only one pebble. Since every outer vertex is of degree $3$, if any vertex is undominated, at least $3$ vertices must be dominated but unpebbled. Hence, in order to obtain an
undominated set of size $\omega+1$, there must be $\omega+4$
vertices that are unpebbled. By the pigeonhole principle, we obtain
a contradiction because there are not enough vertices for this
constraint to hold.  Thus, for $\omega \geq 1$, $n \geq \omega + 3$,
$\Omega_{\omega}(W_n) = n-2-\omega$.
\end{proof}
\section{Graphs of Diameter 2 and 3}

\begin{thm}
Let $G$ be a graph of diameter two with $n$ vertices.  For $\omega
\geq 1$, $\Omega_{\omega}(G) \leq n - 1 - \omega$.

\end{thm}
\begin{proof}
To show that the bound is sharp, consider the graph $H_n$, defined
to be a star graph of order $n$ with $\omega$ additional edges
added to make the graph induced by one subset of $\omega+1$ outer
vertices connected.


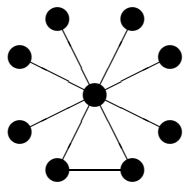
\begin{figure}[htb]
\unitlength 1mm
\begin{center}
\begin{picture}(30,30)
\put(15,15){\circle*{3}} \put(25,20){\circle*{3}}
\put(25,10){\circle*{3}}\put(20,25){\circle*{3}} \put(10,25){\circle*{3}}
\put(5,10){\circle*{3}} \put(20,5){\circle*{3}} \put(10,5){\circle*{3}}
\put(5,20){\circle*{3}}
\put(15,15){\line(1,2){5}} 
\put(15,15){\line(2,1){10}} 
\put(15,15){\line(-2,1){10}} 
\put(15,15){\line(-1,2){5}} 
\put(15,15){\line(1,-2){5}} 
\put(15,15){\line(2,-1){10}} 
\put(15,15){\line(-1,-2){5}} 
\put(15,15){\line(-2,-1){10}} 
\put(10,5){\line(1,0){10}} 

\end{picture}
\end{center}
\caption{An example of the construction for $n = 9$, $\omega = 1$.}
\label{ex2}
\end{figure}


If we place a single pebble on each of the $n-2-\omega$ leaves of
the star that are not connected to any other outer vertices, the
remaining set of undominated vertices is connected and of size
$\omega+1$. Hence, $\Omega(H_n) > n - 2 - \omega$.

Now, let $G$ be a graph of diameter two with $n$ vertices. Suppose
there is an arbitrary configuration of pebbles $c(G)$ that
contains exactly $n-1-\omega$ pebbles. We now show not only that
this configuration  can be solved to eliminate undominated
connected components of order greater than $\omega,$ but can in
fact be solved such that only at most $\omega$ vertices in total
are left undominated.

Much as we did in the proof of Theorem \ref{dia2}, we let
$T_1$ be the set of vertices $v \in G$ such that $c(v)
> 1$, let $T_2$ be the set vertices $w \in G$ such
that $c(w)= 0$ and $w$ is adjacent to some vertex of $T_1,$ and
let $T_3$ be the rest of the vertices, the ones that are neither
in $T_1$ nor adjacent to a vertex of $T_1$. If $|T_3| \leq
\omega$, we are done, because there are no more than $\omega$
undominated vertices and thus the largest undominated component
has size at most $\omega.$ Otherwise, eliminate $\omega$ vertices
in $T_2$ from the graph, and consider the induced subgraph $G'$
and the induced configuration $c'$. We know $G'$ has order $n' = n
- \omega$ and $c'$ still has size at least $n - 1-\omega = n' -1$.
Finally, let $T_1'= T_1,$ $T_2'=T_2$ and $T_3' =T_3 \cap V(G')$.
The new graph $G'$ may no longer have diameter two, which prevents
us from directly applying Theorem \ref{dia2}. Nevertheless, we
notice that in $G',$ every vertex in $T_2'$ is still adjacent to a
vertex in $T_1',$ and every vertex in $T_3'$ is still adjacent to
one in $T_2'.$ Also, since in $G$ we know $d(T_1, T_3)=2,$ it
follows that no path of length one or two between a vertex in
$T_1$ and another vertex of $G$ can pass through $T_3,$ unless
this vertex is the other endpoint. In particular, since the
diameter of $G$ is 2, this implies that the shortest path between
a vertex in $T_1$ and another vertex of $G$ cannot pass through a
vertex of $T_3$ as an intermediate vertex, and so the length of
the shortest path between a vertex in $T_1$ and another vertex in
$G$ will be unaffected by removing a subset of $T_3$. This shows
that in $G',$ if $s \in T_1'$ and $v \in G'$ then $d(s, v) \leq
2$.

We now note that since we have the right number of pebbles in $c'$
(at least $n' - 1$) we can apply the proof of Theorem \ref{dia2}.
Following the proof, we see that we will have $S_1=T_1'$,
$S_2=T_2'$ and $S_3=T_3'.$ Henceforth, the proof never uses the
fact that two vertices of the graph have distance at most two from
one another except when at least one of the vertices in $S_1.$
Thus, the algorithm detailed in the proof can be applied
\emph{mutatis mutandis} to $G'$, after with $G'$ is dominated by
covered vertices. The same sequence of pebbling moves, if
performed on $G,$ leaves all vertices except possibly the $\omega$
that were eliminated to get $G'$ dominated by covered vertices,
thus solving $G$ as desired.
\end{proof}

In general, however, we believe that determining good diameter bounds for $\Omega_w$ will be harder than it is for $\psi.$ It is not even clear to the authors how to construct graphs which establish good lower bounds for large diameters.
However, we conclude this section by conjecturing an analogous result for graphs of
diameter $3$, along with a valid lower-bound construction for this
conjecture.

\begin{con}
Let $G$ be a graph of diameter 3 with $n$ vertices.  For $\omega \geq 1$ ,
$\Omega_\omega(G) \leq \lfloor \frac{3}{2}(n  - 2 - \omega) + 1 \rfloor$.
\end{con}

To see that this result, if true, would give a sharp bound, we exhibit a graph $G$ on $n \geq \omega+3$ vertices such that $\Omega_\omega(G) >
\lfloor \frac{3}{2}(n - 2 - \omega) \rfloor$.   Take
a $K_{\omega + 1}$ and attach  each of its vertices to some other
vertex $v$.  Connect $v$ to each vertex of a $K_{\lceil \frac{n-\omega -
2}{2} \rceil}$, call it $H$. Connect each of the remaining $\lfloor
\frac{n-\omega - 2}{2} \rfloor$ vertices to a vertex of $H$, so that each
vertex in $H$ has at most one such vertex adjacent to it. Now, place three pebbles on
each of the ``tendril" vertices attached to $H$, and if there is one vertex in $H$ without a
tendril, place one pebble on it.  This is a total of $3 \lfloor
\frac{n-\omega - 2}{2} \rfloor$ ($+ 1$ if $n-\omega-2$ is odd) pebbles in
this configuration, which is equivalent to $\lfloor \frac{3}{2}(n  - 2 -
\omega) \rfloor$.  Since it is clearly not possible to dominate the vertices in the
$K_{\omega +1}$, the graph still has an undominated component of order
$\omega + 1$.  Thus, $\Omega_\omega(G) > \lfloor \frac{3}{2}(n  - 2 - \omega)
\rfloor$.


\subsection*{Acknowledgement} The authors would like to thank East Tennessee
State REU director Anant Godbole for supervising this research,
funded by NSF grant DMS-0139286.

\end{document}